\documentclass[12pt]{gITR2e}
\newcommand{\half}{{\textstyle \frac{1}{2}}}
\newcommand{\twothird}{{\textstyle \frac{2}{3}}}
\newcommand{\onethird}{{\textstyle \frac{1}{3}}}
\newcommand{\onesixth}{{\textstyle \frac{1}{6}}}
\newcommand{\threehalves}{{\textstyle \frac{3}{2}}}

\begin{document}
\title{\itshape Hierarchies of sum rules for squares of spherical Bessel functions}
\author{\name{L G Suttorp and A J van  Wonderen}
\affil{Institute for Theoretical Physics, University of Amsterdam,
Science Park 904, 1098 XH Amsterdam, The Netherlands}}
\maketitle
\begin{abstract}
  A four-term recurrence relation for squared spherical Bessel functions is
  shown to yield closed-form expressions for several types of finite
  weighted sums of these functions.  The resulting sum rules, which may
  contain an arbitrarily large number of terms, are found to constitute
  three independent hierarchies. Their use leads to an efficient numerical
  evaluation of these sums.
\end{abstract}
\begin{classcode}
33C10
\end{classcode}

\numberwithin{equation}{section}

\section{Introduction}
Infinite sums of squares of spherical Bessel functions $j_k(z)$ of the form
$\sum_k c_k [j_k(z)]^2$, with various coefficients $c_k$, have been studied
in quite some detail. Several simple examples can be found in standard texts
\cite{AbramowitzStegun68, NIST10}. A more extensive list of such sums has
been compiled in \cite{Luke69}, mainly in the form of expansions of
generalized hypergeometric functions. In contrast, information on finite sums of
this type, with an arbitrary large but finite number of terms, is less well
available. Such finite sums occur in various branches of mathematical
physics, for instance in atomic orbital theory \cite{Harris00}, in acoustic
diffraction problems \cite{Rottbrand00} and more recently in quantum optics
\cite{LGSAJvW15}.

In the present paper we will determine a collection of closed-form
expressions for finite sums of squares of spherical Bessel functions
$j_k(z)$ with $z$-independent coefficients $c_k$. Our main results are
(\ref{12}), (\ref{18}) and (\ref{29}).

\section{Lowest order sum rules} \label{sec2} 
The standard relations connecting three spherical Bessel functions $j_k(z)$
of contiguous order can be employed to derive a four-term recurrence
relation for their squares:
\begin{eqnarray}
   (2k-1)[j_{k-1}(z)]^2-(2k+1)[j_{k}(z)]^2=\rule{5cm}{0mm}\nonumber\\
 =\frac{z^2}{2k-1}\left\{[j_{k-2}(z)]^2-[j_k(z)]^2\right\}
+\frac{z^2}{2k+1}\left\{[j_{k-1}(z)]^2-[j_{k+1}(z)]^2\right\}\,
 . \label{1}
\end{eqnarray}
In fact, the proof follows by elimination of $j_{k-2}(z)$ and $j_{k+1}(z)$
in favour of $j_{k-1}(z)$ and $j_k(z)$ (see \cite{AbramowitzStegun68},
formula 10.1.19). Multiplying this identity by an as yet undetermined
coefficient $a_k$, summing over $k$ and shifting the summation variables,
one finds a relation between finite sums:
\begin{eqnarray}
   \sum_{k=0}^\ell f_k[j_{k}(z)]^2=z^2 \sum_{k=0}^\ell g_k [j_{k}(z)]^2
 -z^2\frac{1}{2\ell+3}(a_{\ell+1}+a_{\ell+2})[j_{\ell}(z)]^2\nonumber\\
-z^2\frac{1}{2\ell+1}(a_{\ell+1}+a_\ell)[j_{\ell+1}(z)]^2
+2z\, a_{\ell+1}\, j_{\ell}(z)j_{\ell+1}(z)+F(z)\, , \label{2}
\end{eqnarray}
with $\ell \geq 0$, $f_k=(2k+1)(a_{k+1}-a_k)$, and $g_k=
(a_{k+2}+a_{k+1})/(2k+3)-(a_k+a_{k-1})/(2k-1)$. The $\ell$-dependent terms
at the right-hand side follow upon adjusting the upper bounds of the
sums. The last term results by completing the sums at their lower bounds;
it can easily be found by putting $\ell=0$.  The term with $j_\ell(z)
j_{\ell+1}(z)$ arises upon using the equality
\begin{eqnarray}
   z j_\ell(z) j_{\ell+1}(z)=\half
z^2\frac{1}{2\ell+1}\left\{[j_{\ell+1}(z)]^2-[j_{\ell-1}(z)]^2\right\}
+\half(2\ell+1)[j_\ell(z)]^2\, , \label{3}
\end{eqnarray}
which is established by elimination of $j_{\ell-1}(z)$ in a similar way as in
(\ref{1}).

The identity (\ref{2}) gives a relation between two sums. In general, it
does not yield a closed-form expression for any of these. However, one may
derive explicit sum rules in two different ways: either by choosing $a_k$
in such a way that the left-hand side vanishes, or such that the sum at the
right-hand side drops out.  In the first case the coefficients $a_k$
should satisfy the relation $f_k=0$ or $a_{k+1}=a_k$, so that one may take
$a_k=1$ for all $k$. In this way one arrives at the sum rule:
\begin{eqnarray}
   \sum_{k=0}^\ell \frac{1}{(2k-1)(2k+3)} [j_k(z)]^2=
-\frac{1}{4(2\ell+3)}[j_\ell(z)]^2 -\frac{1}{4(2\ell+1)}[j_{\ell+1}(z)]^2\nonumber\\
+\frac{1}{4z}j_\ell(z)j_{\ell+1}(z)-\frac{1}{2z}j_1(2z)\, , \label{4}
\end{eqnarray}
where the last term is found by choosing $\ell=0$.

In the other case the coefficients $a_k$ must fulfill the condition
$g_k=0$ for all $k$. Solving the ensuing recurrence relation for $a_k$, by
first introducing $b_k=a_{k-1}+a_k$ and evaluating $b_k$ for even and odd
$k$ separately, one finds $a_k$ in terms of the initial conditions $a_0$,
$a_1$ and $a_2$ as
\begin{equation}
   a_k=[\half k+(-1)^{k+1}(\half k^2-1)]a_0+[\twothird k+\onethird(-1)^{k+1} k^2]a_1
+[\onesixth k+\onesixth (-1)^k k^2]a_2 \, . \label{5}
\end{equation}
As a consequence, the coefficient $f_k$ in the sum at the left-hand side of (2)
gets the form
\begin{eqnarray}
   f_k=(2k+1)\left\{[\half+(-1)^k(k^2+k-\threehalves)]a_0+
[\twothird+(-1)^k(\twothird k^2+\twothird k+\onethird)]a_1\right.\nonumber\\
\left.+[\onesixth+(-1)^{k+1}(\onethird k^2+\onethird k+\onesixth)]a_2\right\}\, . \label{6}
\end{eqnarray}

By making specific choices for $a_0$, $a_1$ and $a_2$ one may arrive at sum
rules with either alternating coefficients (proportional to $(-1)^k$) or
non-alternating coefficients. The latter type shows up by choosing $a_0=0$
and $a_2=2a_1$. Taking $a_1=1$ one gets $a_k=k$ and $a_{k+1}-a_k=1$. As a
consequence, one arrives at a second sum rule with non-alternating
coefficients:
\begin{equation}
   \sum_{k=0}^\ell (2k+1)[j_k(z)]^2=-z^2 [j_\ell(z)]^2-z^2
[j_{\ell+1}(z)]^2+2(\ell+1)z j_\ell(z) j_{\ell+1}(z) +1 \, , \label{7}
\end{equation} 
where the last term is obtained by taking $\ell=0$, as before. 

Further sum rules, with alternating coefficients, arise by choosing
$a_2=-3a_0-4a_1$, while $a_0$ and $a_1$ can still be chosen
at will. One possible choice is $a_0=-a_1=1$, which implies $a_k=(-1)^k$. The
ensuing sum rule following from (\ref{2}) is
\begin{equation}
\sum_{k=0}^\ell (-1)^k(2k+1)[j_k(z)]^2=(-1)^\ell z
j_\ell(z)j_{\ell+1}(z)+j_0(2z)\, . \label{8}
\end{equation}

Finally, one may take $a_k$ to be proportional to a quadratic polynomial in
$k$. Upon putting $a_0=a_1=1$ we get $a_k=(-1)^{k+1} (2k^2-1)$. This choice
yields the last sum rule of our set:
\begin{eqnarray}
   \sum_{k=0}^\ell (-1)^k(2k+1)k(k+1)[j_k(z)]^2=\half (-1)^\ell z^2
[j_\ell(z)]^2-\half (-1)^\ell z^2 [j_{\ell+1}(z)]^2\nonumber\\
+(-1)^\ell (\ell^2+2\ell+\half) z
j_\ell(z) j_{\ell+1}(z)-z j_1(2z) \, . \label{9}
\end{eqnarray}

In conclusion, by making judicious choices for the coefficients $a_k$ in
(\ref{2}) we have derived several independent sum rules for the squares of
the spherical Bessel functions: two rules with non-alternating
coefficients, namely (\ref{4}) and (\ref{7}), and two closely related ones
with alternating coefficients, namely (\ref{8}) and (\ref{9}). In the
following we shall see that these sum rules can be used as a basis from
which three independent hierarchies of sum rules can be established.

\section{Hierarchies of sum rules}
The sum rules (\ref{4}) and (\ref{7})--(\ref{9}) are part of several
hierarchies of sum rules. These hierarchies follow by using (\ref{2}) for
suitable $a_k$ as a recurrence relation. As a first example we shall start
from (\ref{4}), by choosing the coefficient $a_k$ in (\ref{2}) such that at
its left-hand side the sum found in (\ref{4}) shows up. Apart from a
trivial factor this implies that $a_k$ should fulfill the relation
$a_{k+1}-a_k=1/(k-\half)_3$, with $(f)_n=f(f+1)\ldots(f+n-1)$ the
Pochhammer symbol. Solving for $a_k$ we find $a_k=-1/[2(k-\half)_2]$, where
we chose the initial condition as $a_0=2$.  Inserting this form for $a_k$
in the sum at the right-hand side of (\ref{2}) we obtain as its coefficient
$g_k=3(2k+1)/[2(k-\threehalves)_5]$.  Hence, the relation (\ref{2}) yields
an expression for a sum with a coefficient proportional to
$(2k+1)/(k-\threehalves)_5$, which is the next in a hierarchy of sum rules
of which (\ref{4}) was the first. In fact, the above procedure can be
repeated. By choosing $a_k^{(p)}$ in (\ref{2}) in such a way that the
coefficient in the sum at the left-hand side is
$f_k^{(p)}=(2k+1)/(k-p-\half)_{2p+3}$ (for arbitrary integer $p\geq 0$), we
arrive at a coefficient $g_k^{(p)}$ at the right-hand side that is
proportional to $(2k+1)/(k-p-\threehalves)_{2p+5}$, for a suitable choice of
the initial condition. To achieve this one should take
\begin{equation}
a_k^{(p)}=-\frac{1}{2(p+1)\, (k-p-\half)_{2p+2}} \, , \label{10}
\end{equation} 
for $p\geq 0$ and all $k$. With this choice of $a_k^{(p)}$ the relation
(\ref{2}) becomes for $p\geq 0$ and $\ell\geq 0$:
\begin{eqnarray}
   \sum_{k=0}^\ell \frac{2k+1}{(k-p-\threehalves)_{2p+5}} [j_k(z)]^2=
\frac{2(p+1)}{z^2(2p+3)}
\sum_{k=0}^\ell\frac{2k+1}{(k-p-\half)_{2p+3}} [j_k(z)]^2 \rule{2cm}{0mm}\nonumber\\
   -\frac{1}{(2p+3)\, (\ell-p+\half)_{2p+3}}[j_\ell(z)]^2
-\frac{1}{(2p+3)\, (\ell-p-\half)_{2p+3}}[j_{\ell+1}(z)]^2\rule{2cm}{0mm}\nonumber\\
   +\frac{2}{z(2p+3)\, (\ell-p+\half)_{2p+2}}j_\ell(z)j_{\ell+1}(z)
+\frac{1}{(-p-\threehalves)_{2p+4}}
\left(\frac{p+1}{z^2}  j_0(2z)+\frac{1}{z}j_1(2z)\right) , \quad
\label{11}
\end{eqnarray}
where the last term has been determined by putting $\ell=0$. Upon using
this identity repeatedly we arrive at a {\em first} hierarchy of sum rules
of the form:
\begin{eqnarray}
   \sum_{k=0}^\ell \frac{2k+1}{(k-p-\half)_{2p+3}}[j_k(z)]^2=
z^2\, A_\ell^{[1],(p)}(z)\, [j_\ell(z)]^2+
z^2\, B_\ell^{[1],(p)}(z)\, [j_{\ell+1}(z)]^2\nonumber\\
   +  z\, C_\ell^{[1],(p)}(z)\, j_\ell(z) j_{\ell+1}(z)
+\frac{1}{(-p-\half)_{2p+2}}\sum_{k=0}^{p}
(-1)^{k}\frac{(p-k+1)_{k}}{z^{k+1}} j_{k+1}(2z)\, , \label{12}
\end{eqnarray}
for any $p\geq 0$ and $\ell \geq 0$. The coefficients at the right-hand
side are polynomials in $1/z^2$:
\begin{eqnarray}
   A_\ell^{[1],(p)}(z)=B_{\ell+1}^{[1],(p)}(z)=
-\half\sum_{k=0}^p\frac{(p-k+1)_k}{(p-k+\half)_{k+1}(\ell-p+k+\threehalves)_{2p-2k+1}}
\frac{1}{z^{2k+2}}\, , \label{13}\\
   C_\ell^{[1],(p)}(z)=\sum_{k=0}^p\frac{(p-k+1)_k}{(p-k+\half)_{k+1}(\ell-p+k+
\threehalves)_{2p-2k}}\frac{1}{z^{2k+2}}\, . \label{14}
\end{eqnarray}
The last term in (\ref{12}) is obtained from (\ref{11}) by using the
recurrence relations for $j_k(2z)$.  For $p=0$ one recovers the sum rule
(\ref{4}), which we have used as our starting-point. For small $p>0$ and
arbitrary $\ell$ the coefficients of the squared spherical Bessel
functions at the right-hand side of (\ref{12}) are small-degree polynomials
in $1/z$ that are easily evaluated. For general $p\geq 0$ the sum rules
(\ref{12}), with spherical Hankel functions instead of spherical Bessel
functions, were of crucial importance in the analysis of the modified
atomic decay rates in \cite{LGSAJvW15}.

A rather different hierarchy follows by starting from the sum rule (\ref{7}) and
choosing $g_k=2k+1$ in (\ref{2}). Solving for $a_k$ one finds $a_k=\half
(k-1)_3$ for all $k$, when a suitable choice of initial conditions is
made.  Subsequently, $f_k$ is obtained as $f_k=\threehalves (2k+1) (k)_2$,
so that (\ref{2}), with (\ref{7}) inserted at the right-hand side, leads to
a sum rule for $\ell$-dependent sums with a coefficient $(2k+1)(k)_2$. The
procedure can be generalized by taking
\begin{equation}
a_k^{(p)}=\frac{(k-p-1)_{2p+3}}{2(p+1)} \, , \label{15}
\end{equation}
for $p\geq 0$ and all $k$, and hence
\begin{equation}
 f_k^{(p)}=\frac{2p+3}{2(p+1)}(2k+1)(k-p)_{2p+2} \quad , \quad 
g_k^{(p)}=(2k+1)(k-p+1)_{2p}\, . \label{16}
\end{equation}
In this way we get from (\ref{2}) a relation between sums of similar form:
\begin{eqnarray}
   \frac{2p+3}{2(p+1)}\sum_{k=p+1}^\ell (2k+1)(k-p)_{2p+2}[j_k(z)]^2=
z^2\sum_{k=p}^\ell (2k+1)(k-p+1)_{2p}[j_k(z)]^2\rule{2cm}{0mm}\nonumber\\
-\frac{z^2}{2(p+1)}(\ell-p+1)_{2p+2}[j_\ell(z)]^2-
\frac{z^2}{2(p+1)}(\ell-p)_{2p+2}[j_{\ell+1}(z)]^2\rule{3cm}{0mm}\nonumber\\
+\frac{z}{p+1}(\ell-p)_{2p+3}\, j_\ell(z)j_{\ell+1}(z)\, , \rule{3cm}{0mm}\label{17}
\end{eqnarray}
for $\ell\geq p \geq 0$. The last term in (\ref{2}) is found to be $0$ in this case,
as follows by putting $\ell=p$. By employing this identity recursively and
using (\ref{7}), one arrives at a {\em second} hierarchy of sum rules with
non-alternating coefficients, on a par with (\ref{12}):
\begin{eqnarray}
   \sum_{k=p}^\ell (2k+1)(k-p+1)_{2p}[j_k(z)]^2=
z^2\, A_\ell^{[2],(p)}(z)\, [j_\ell(z)]^2+
z^2\, B_\ell^{[2],(p)}(z)\, [j_{\ell+1}(z)]^2\nonumber\\
+  z\, C_\ell^{[2],(p)}(z)\, j_\ell(z) j_{\ell+1}(z)
+\frac{p!}{(\threehalves)_p}\, z^{2p}\, , \label{18}
\end{eqnarray}
for all $\ell\geq p\geq 0$. The coefficients at the right-hand side
are polynomials in $z^2$:
\begin{eqnarray}
A_{\ell}^{[2],(p)}(z)=B_{\ell+1}^{[2],(p)}(z)=-\half\sum_{k=0}^p \frac{(p-k+1)_k\, (\ell
    -p+k+2)_{2p-2k}}{(p-k+\half)_{k+1}} z^{2k}\, ,\label{19}\\
C_{\ell}^{[2],(p)}(z)=\sum_{k=0}^p \frac{(p-k+1)_k\, (\ell
    -p+k+1)_{2p-2k+1}}{(p-k+\half)_{k+1}} z^{2k}\, . \label{20}
\end{eqnarray}
For $p=0$ the sum rule (\ref{18}) reduces to (\ref{7}), which served as the
basis of the hierarchy. For small $p>0$ the $z$-dependent polynomials
occurring at the right-hand side of (\ref{18}) have got a small degree, as in
(\ref{12}). Sum rules closely related to (\ref{18}), with spherical Hankel
functions as before, have been used in the quantum optics problem in
\cite{LGSAJvW15}. Comparing the two hierarchies (\ref{12}) and (\ref{18})
we see that the coefficients in the weighted sums differ
considerably. Whereas the coefficient in (\ref{12}) contains an odd number
of factors in the denominator, the Pochhammer symbol in (\ref{18}) leads to
an even number of factors in the numerator, at least if the common factor
$2k+1$ is left out of consideration. Furthermore, the factors in (\ref{12})
are half-integer (so that no singularities can arise), and in (\ref{18})
they are all integer. In both cases the number of factors increases with
$p$.

Finally, we can build a hierarchy of sum rules with alternating
coefficients by starting from (\ref{8}) and (\ref{9}). When we choose $a_k$ to
have the somewhat elaborate form
\begin{equation}
a_k^{(p)}=(-1)^k \frac{k}{2}\sum_{m=0}^p c_m^{(p)} \frac{1}{(m+1)(m+2)}(k-m-1)_{2m+3} \, , \label{21}
\end{equation}
with 
\begin{equation}
c_m^{(p)}=\frac{1}{2^{2p-2m}m!(p-m)!} (2m-p+2)_{2p-2m} \, , \label{22}
\end{equation}
the coefficients $f_k$ and $g_k$ are found as
\begin{eqnarray}
f_k^{(p)}=(-1)^{k+1}(2k+1)\sum_{m=0}^{p+2} c_m^{(p+2)} \, (k-m+1)_{2m} \, , \label{23}
\\
g_k^{(p)}=(-1)^k(2k+1)\sum_{m=0}^p c_m^{(p)}\, (k-m+1)_{2m}\, , \label{24}
\end{eqnarray}
for all $p\geq 0$ and $k\geq 0$. In deriving (\ref{23}) we have used the
recurrence relation $c_m^{(p+2)}=c_{m-2}^{(p)}/[m(m-1)]+c_{m-1}^{(p)}\,
(2m+1)/(2m)$ for $m\geq 2$ and $p\geq 0$.  Since the expressions (\ref{23})
and (\ref{24}) are closely analogous, with $f_k^{(p)}=-g_k^{(p+2)}$, one
may use (\ref{2}) to derive a recurrence relation for sums of a similar
type:
\begin{eqnarray}
\sum_{k=0}^\ell (-1)^k (2k+1)\left[\sum_{m=0}^{p+2} c_m^{(p+2)}\,
  (k-m+1)_{2m}\right] [j_k(z)]^2=\nonumber\\
 = -z^2\sum_{k=0}^\ell (-1)^k (2k+1)\left[\sum_{m=0}^{p} c_m^{(p)}\,
  (k-m+1)_{2m}\right] [j_k(z)]^2\nonumber\\
 +\half (-1)^\ell  z^2\left[\sum_{m=0}^{p} c_m^{(p)}\,\frac{1}{m+1}
  (\ell-m+1)_{2m+2}\right] [j_\ell(z)]^2\nonumber\\
 -\half (-1)^\ell  z^2\left[\sum_{m=0}^{p} c_m^{(p)}\,\frac{1}{m+1}
  (\ell-m)_{2m+2}\right] [j_{\ell+1}(z)]^2\nonumber\\
 +(-1)^\ell (\ell+1)\left[\sum_{m=0}^{p} c_m^{(p)}\,\frac{1}{(m+1)(m+2)}
  (\ell-m)_{2m+3}\right] z j_\ell(z) j_{\ell+1}(z) \, . \label{25}
\end{eqnarray}
Once again the last term in (\ref{2}) drops out, as follows by taking
$\ell=0$. Starting from the sum rules (\ref{8}) and (\ref{9}) and using
the recurrence relation separately for even and odd values of $p$ we may
obtain explicit sum rules for all $p$. For any $p\geq 0$ and $\ell\geq0$ we
get:
\begin{eqnarray}
   \sum_{k=0}^\ell (-1)^k (2k+1)\left[\sum_{m=0}^{p} c_m^{(p)}\,
  (k-m+1)_{2m}\right] [j_k(z)]^2=\rule{5cm}{0mm}\nonumber\\
    =\half (-1)^\ell \left[ \sum_{m=1}^{[p/2]} (-1)^{m+1} z^{2m}
\sum_{n=0}^{p-2m} c_n^{(p-2m)}\frac{1}{n+1}(\ell-n+1)_{2n+2}
+\delta_p^o \, (-1)^{(p-1)/2} z^{p+1}\right] 
[j_\ell(z)]^2\nonumber\\
    +\half (-1)^\ell \left[ \sum_{m=1}^{[p/2]} (-1)^m z^{2m}
\sum_{n=0}^{p-2m} c_n^{(p-2m)}\frac{1}{n+1}(\ell-n)_{2n+2}
+\delta_p^o\, (-1)^{(p+1)/2} z^{p+1}\right] 
[j_{\ell+1}(z)]^2\rule{7mm}{0mm}\nonumber\\
    + (-1)^\ell \left[ (\ell+1) \sum_{m=1}^{[p/2]} (-1)^{m+1} z^{2m-1}
\sum_{n=0}^{p-2m}
c_n^{(p-2m)}\frac{1}{(n+1)(n+2)}(\ell-n)_{2n+3}\right.\rule{2cm}{0mm}\nonumber\\
 \left.+\delta_p^e \, (-1)^{p/2} z^{p+1}+\delta_p^o \, (-1)^{(p-1)/2} z^p (\ell+1)^2\rule{0mm}{8mm}\right] 
j_\ell(z)j_{\ell+1}(z) \rule{2cm}{0mm}\nonumber\\
   +\delta_p^e \, (-1)^{p/2}z^p j_0(2z) 
+\delta_p^o \, (-1)^{(p-1)/2} z^{p-1} [\half j_0(2z)-z j_1(2z)]\, . \rule{1cm}{0mm}\label{26}
\end{eqnarray}
Here $\delta_p^e$ equals 1 for even $p$, and 0 for odd $p$, while
$\delta_p^o$ is defined analogously, with even and odd interchanged. The
upper bounds of the summations contain the `entier' function $[x]$ which is
the largest integer $\leq x$.  For $p=0$ the sum rule (\ref{26}) reduces to
(\ref{8}), whereas for $p=1$ a linear combination of (\ref{8}) and
(\ref{9}) is recovered.

By taking suitable linear combinations of the sum rules (\ref{26}) for
various values of $p$ we may obtain expressions for sums with the simple
coefficients $(-1)^k (2k+1)(k-p+1)_{2p}$. In fact, we may use the identity
for $q\geq m$:
\begin{equation}
\sum_{p=m}^q f_p^{(q)}\, c_m^{(p)} =\delta_{m,q} \, \quad {\rm with} \quad 
f_p^{(q)}=(-1)^{p+q} \frac{q! (2q-p)!}{2^{2q-2p} p!(q-p)!}\, , \label{27}
\end{equation}
which can be proved by employing a relation (due to Dzhrbashyan
\cite{Dzhrbashyan64}) for a terminating generalized hypergeometric function
${}_3F_2(1)$ with unit argument (see also \cite{Luke69}, formula
3.13.3(9)). We now take the sum $\sum_{p=0}^q f_p^{(q)}$ of (\ref{26}) and
use (\ref{27}) at the left-hand side. In the first term at the right-hand
side we interchange the order of the summations in such a way that the sum over
$p$ can be carried out first. Substitution of (\ref{22}) then leads to an
expression that can be evaluated with the help of \cite{Dzhrbashyan64}:
\begin{eqnarray}
   \sum_{p=2m+n}^q (-1)^p\, \frac{(2q-p)!}{p!(q-p)!(p-2m-n)!}\,
(2n-p+2m+2)_{2p-4m-2n}= \nonumber\\
   =(-1)^n\, \frac{(2q-2m-n)!}{(2m+n)!(q-2m-n)!}\,  \rule{5cm}{0mm}\nonumber\\
\times \mbox{}_3F_2(-q+2m+n,n+2,-n-1;-2q+2m+n,2m+n+1;1)=\nonumber\\
   = (-1)^n\, \frac{2^{2q-4m-2n}}{(m-1)!}\, (q-2m-n+1)_{m-1}\,
(m+n+\threehalves)_{q-2m-n}\, . \label{28}
\end{eqnarray}
Treating the second and third terms in (\ref{26}) in a similar way and
relabelling $q$ as $p$ we find a {\em third} hierarchy of sum rules:
\begin{eqnarray}
   \sum_{k=p}^\ell (-1)^k (2k+1)(k-p+1)_{2p} \, [j_k(z)]^2=
z^2\, A_\ell^{[3],(p)}(z)\, [j_\ell(z)]^2+
z^2\, B_\ell^{[3],(p)}(z)\, [j_{\ell+1}(z)]^2\nonumber\\
+  z\, C_\ell^{[3],(p)}(z)\, j_\ell(z) j_{\ell+1}(z)+
(-1)^p p! \, z^p j_p(2z) \, , \rule{1cm}{0mm}\label{29}
\end{eqnarray}
for $p\geq 0$ and $\ell\geq p$.  The coefficients at the right-hand side
are polynomials in $z^2$ that are given as follows:
\begin{eqnarray}
   A_\ell^{[3],(p)}(z)=B_{\ell+1}^{[3],(p)}(z)=
\half (-1)^{p+\ell+1} p!\sum_{m=0}^{[(p-1)/2]} \sum_{n=0}^{p-2m-1} 
(-1)^{m+n} \frac{1}{m!\, n!}  \nonumber\\
\times (m+n+\threehalves)_{p-2m-n-1} \, (p-2m-n)_m \, (\ell-n+2)_{2n} \, z^{2m}\, , \label{30}\\
   C_\ell^{[3],(p)}(z)= (-1)^{p+\ell} p!\, (\ell+1)\sum_{m=0}^{[p/2]} \sum_{n=0}^{p-2m}
(-1)^{m+n}\frac{1}{m!\, n!} \nonumber\\
\times (m+n+\half)_{p-2m-n}\, (p-2m-n+1)_m\,(\ell-n+2)_{2n-1} \, z^{2m} \, , \label{31}
\end{eqnarray}
with $(f)_{-1}=\Gamma(f-1)/\Gamma(f)=1/(f-1)$ for $f\neq 1$. For $p=0$ and
$p=1$ the sum rule (\ref{29}) yields (\ref{8}) and (\ref{9}),
respectively. The $z$-dependent polynomials in (\ref{29}) are easily
evaluated for small $p>1$ and arbitrary $\ell$, as their degree is small in
that case.

The three hierarchies of sum rules (\ref{12}), (\ref{18}) and (\ref{29}) for
squares of spherical Bessel functions are the main results of this
paper. They have been derived in a systematic way from the four-term
recurrence relation (\ref{1}). 

\section{Discussion and conclusion}

The derivation of the hierarchies for finite sums of squares of spherical
Bessel functions shows that these hierarchies appear to be uniquely defined
as generalizations of the lowest order sum rules from Section
\ref{sec2}. The latter followed from the fundamental recurrence relation
(\ref{2}) for squares of spherical Bessel functions. 

The finite sums of Bessel functions found above converge as the upper limit
$\ell$ tends to $\infty$, since $[j_k(z)]^2$ goes to $0$ quite fast for
$k\rightarrow \infty$ at fixed $z$. For infinite $\ell$ the sum rules (\ref{12}),
(\ref{18}) and (\ref{29}) are consistent with those found before (see
\cite{Luke69}, formulas 9.4.4(13) and 9.4.7(9)). For finite $\ell$ the
sum rules derived above are all new, to the best of our knowledge.

Hierarchies of sum rules similar to those given in (\ref{12}), (\ref{18})
and (\ref{29}) may be established for sums over products
$f_\ell(z)g_\ell(z)$, with $f_\ell(z)$ and $g_\ell(z)$ equal to
$j_\ell(z)$, $y_\ell(z)$, $h_\ell^{(1)}(z)$ or $h_\ell^{(2)}(z)$,
independently, since the recurrence relation (\ref{1}) holds true for any
product of these functions. It should be noted that the term $F(z)$ in
(\ref{2}), and hence the terms independent of $\ell$ in (\ref{12}),
(\ref{18}) and (\ref{29}) get a different form upon switching to general
functions $f_\ell(z)$ and $g_\ell(z)$.  Furthermore, the product $j_\ell(z)
j_{\ell+1}(z)$ must be replaced by $\half [f_\ell(z)
g_{\ell+1}(z)+f_{\ell+1}(z) g_\ell(z)]$. In \cite{LGSAJvW15} the sum rules
(\ref{12}) and (\ref{18}), with spherical Hankel functions instead of
$j_\ell(z)$, have been used to determine indefinite integrals over squares
of these functions. The general form of these indefinite integrals is
$\int^z du\, u^{-n}\, h_{\ell_1}^{(1)}(u)\,h_{\ell_2}^{(i)}(u)$, for real
$z>0$, integer $n$, non-negative integers $\ell_1,\ell_2$ and
$i=1,2$. 

The identities (\ref{12}), (\ref{18}) and (\ref{29}) yield an efficient way
to evaluate the sums at their left-hand sides numerically, in particular
for small $p$ and large $\ell$. Comparing for instance the evaluation times
of both sides of (\ref{29}) for $p=0$, $\ell=50$ and $z=50$ with the help
of the numerical software contained in Mathematica, one finds that
calculating the right-hand side is more than 10 times faster than
calculating the left-hand side. Such an increase in the efficiency of the
numerical evaluation proved to be advantageous in producing the plots in
\cite{LGSAJvW15} for a range of values of $z$.

Remarkably enough, the derivation presented above shows the close
connection between the three hierarchies (\ref{12}), (\ref{18}) and
(\ref{29}): all three follow, on an equal footing, from the single fundamental
recurrence relation (\ref{1})  for the squared spherical Bessel functions.

\end{document}